\def\yes{\if00}
\def\no{\if01}
\def\iftwelvept{\yes}
\def\ifusepdf{\no}
\def\ifpsfont{\yes}

\iftwelvept
\documentclass[leqno,12pt]{amsart}
\else
\documentclass[leqno]{amsart}
\fi
\usepackage{amssymb}
\usepackage{amscd}
\usepackage{latexsym}
\usepackage{verbatim}
\ifusepdf
\usepackage{hyperref}
\else\fi
\ifpsfont
\usepackage[T1]{fontenc}
\usepackage{times}
\else\fi

\iftwelvept
\else\fi

\theoremstyle{plain}
\newtheorem{Theorem}{Theorem}[section]
\newtheorem{Proposition}[Theorem]{Proposition}
\newtheorem{Lemma}[Theorem]{Lemma}
\newtheorem{Corollary}[Theorem]{Corollary}

\newtheorem{Claim}{Claim}[Theorem]

\theoremstyle{definition}

\newtheorem{Remark}[Theorem]{Remark}

\renewcommand{\theTheorem}{\arabic{section}.\arabic{Theorem}}
\renewcommand{\theClaim}{\arabic{section}.\arabic{Theorem}.\arabic{Claim}}
\renewcommand{\theequation}{\arabic{section}.\arabic{Theorem}.\arabic{Claim}}

\def\rom{\textup}
\newcommand{\ZZ}{{\mathbb{Z}}}
\newcommand{\QQ}{{\mathbb{Q}}}
\newcommand{\RR}{{\mathbb{R}}}
\newcommand{\CC}{{\mathbb{C}}}

\newcommand{\PP}{{\mathbb{P}}}

\newcommand{\OO}{{\mathcal{O}}}
\newcommand{\XX}{{\mathcal{X}}}
\newcommand{\YY}{{\mathcal{Y}}}
\newcommand{\ZZZ}{{\mathcal{Z}}}
\newcommand{\LL}{{\mathcal{L}}}
\newcommand{\MM}{{\mathcal{M}}}
\newcommand{\codim}{\operatorname{codim}}
\newcommand{\ord}{\operatorname{ord}}

\newcommand{\Pic}{\operatorname{Pic}}
\newcommand{\Gal}{\operatorname{Gal}}

\newcommand{\Spec}{\operatorname{Spec}}
\newcommand{\Chow}{\operatorname{Ch}}

\newcommand{\Proj}{\operatorname{Proj}}
\newcommand{\zero}{\operatorname{div}}
\newcommand{\Proof}{{\sl Proof.}\quad}
\newcommand{\adeg}{\widehat{\operatorname{deg}}}
\newcommand{\trdeg}{\operatorname{tr.deg}}
\newcommand{\rank}{\operatorname{rk}}
\newcommand{\acherncl}{\widehat{{c}}}

\newcommand{\QED}{{\unskip\nobreak\hfil\penalty50\quad\null\nobreak\hfil
{$\Box$}\parfillskip0pt\finalhyphendemerits0\par\medskip}}
\newcommand{\rest}[2]{\left.{#1}\right\vert_{{#2}}}

\begin{document}

\title[A note on polarizations of finitely generated fields]%
{A note on polarizations of finitely generated fields}
\author{Atsushi Moriwaki}
\address{Department of Mathematics, Faculty of Science,
Kyoto University, Kyoto, 606-8502, Japan}
\email{moriwaki@kusm.kyoto-u.ac.jp}
\date{4/June/2000, 10:30PM (JP), (Version 1.0)}


\maketitle

\renewcommand{\theTheorem}{\Roman{Theorem}}
\section*{Introduction}
In the paper \cite{MoArht}, we established Northcott's theorem for
height functions over finitely generated fields.
Unfortunately, Northcott's theorem on finitely generated fields
does not hold in general (cf. Remark~\ref{rem:north:not:hold}).
Actually, it depends on the choice of a polarization.
In this short note, we will propose a weaker condition of the polarization to
guarantee Northcott's theorem.
We will also show the generalization of conjectures of Bogomolov and Lang 
in \cite{MoBL} under the weaker polarization.

First of all, let us introduce the weaker condition of
a polarization.
Let $K$ be a finitely generated field over $\QQ$.
A polarization $\overline{B} = (B; \overline{H}_1, \ldots, \overline{H}_d)$
of $K$
is said to be {\em fairly large}
if there are generically finite morphisms
$\mu : B' \to B$ and $\nu : B' \to \left(\PP^1_{\ZZ}\right)^d$ of
flat and projective integral schemes over $\ZZ$, and
nef and big $C^{\infty}$-hermitian $\QQ$-line bundles
$\overline{L}_1, \ldots, \overline{L}_d$
on $\PP^1_{\ZZ}$
such that a positive power of
$\mu^*(\overline{H}_i) \otimes 
\nu^*(p_i^*(\overline{L}_i))^{\otimes -1}$ has a small section for every $i$,
where $p_i : \left(\PP^1_{\ZZ}\right)^d \to \PP^1_{\ZZ}$ is the projection
to the $i$-th factor.
Then, we have the following Theorem~\ref{thm:northcott:intro},
Theorem~\ref{thm:BL:conj:fair:pol:intro} and
Theorem~\ref{thm:bogo:equiv:height},
which are generalizations of the previous results.

\begin{Theorem}[{\cite[Theorem~4.3]{MoArht}}]
\label{thm:northcott:intro}
We assume that the polarization
$\overline{B}$ is fairly large.
Let $X$ be a geometrically irreducible projective variety over $K$, and
$L$ an ample line bundle on $X$.
Then, for any number $M$ and any positive integer $e$, the set
\[
\{ x \in X(\overline{K}) \mid h^{\overline{B}}_L(x) \leq M, \quad
[K(x) : K] \leq e \}
\]
is finite.
Moreover, it can be
generalized to the height of cycles on a projective variety
\rom{(}cf. \rom{Theorem~\ref{thm:northcott:cycle:intro}}\rom{)}.
\end{Theorem}

\begin{Theorem}[{\cite[Theorem~A]{MoBL}}]
\label{thm:BL:conj:fair:pol:intro}
We assume that the polarization
$\overline{B}$ is fairly large.
Let $A$ be an abelian variety over $K$, and $L$ a symmetric
ample line bundle on $A$.
Let 
\[
\langle\ , \ \rangle_L^{\overline{B}} :
A(\overline{K}) \times A(\overline{K}) \to \RR
\]
be a paring given by
\[
\langle x , y \rangle_L^{\overline{B}} =
\frac{1}{2} \left( \hat{h}_L^{\overline{B}}(x+y)
-  \hat{h}_L^{\overline{B}}(x) -  \hat{h}_L^{\overline{B}}(x) \right).
\]
For $x_1, \ldots, x_l \in A(\overline{K})$,
we denote $\det \left( \langle x_i, x_j \rangle_L^{\overline{B}} \right)$
by $\delta_L^{\overline{B}}(x_1, \ldots, x_l)$.

Let $\Gamma$ be a subgroup of finite rank in $A(\overline{K})$, and
$X$ a subvariety of $A_{\overline{K}}$.
Fix a basis $\{\gamma_1, \ldots, \gamma_n \}$ of $\Gamma \otimes \QQ$.
If the set $\{ x \in X(\overline{K}) \mid
\delta_L^{\overline{B}}(\gamma_1, \ldots, \gamma_n, x) \leq \epsilon \}$
is Zariski dense in $X$ for every positive number $\epsilon$,
then $X$ is a translation of an abelian subvariety of $A_{\overline{K}}$
by an element of $\Gamma_{div}$,
where $\Gamma_{div} = \{ x \in A(\overline{K}) \mid
\text{$nx \in \Gamma$ for some positive integer $n$} \}$.
\end{Theorem}

\begin{Theorem}[{\cite[Theorem~5.1]{MoCanht}}]
\label{thm:bogo:equiv:height}
We assume that the polarization
$\overline{B}$ is fairly large.
For a subvariety $X$ of $A_{\overline{K}}$,
the following are equivalent.
\begin{enumerate}
\renewcommand{\labelenumi}{(\arabic{enumi})}
\item
$X$ is a translation of an abelian subvariety by a torsion point.

\item
The set 
$\{ x \in X(\overline{K}) \mid \hat{h}_L^{\overline{B}}(x)
   \leq \epsilon \}$
is Zariski dense in $X$ for every $\epsilon > 0$.

\item
The canonical height of $X$ with respect to $L$ and $\overline{B}$
is zero, i.e.,
$\hat{h}^{\overline{B}}_{L}(X) = 0$.
\end{enumerate}
\end{Theorem}

\renewcommand{\theTheorem}{\arabic{section}.\arabic{Theorem}}
\renewcommand{\theClaim}{\arabic{section}.\arabic{Theorem}.\arabic{Claim}}
\renewcommand{\theequation}{\arabic{section}.\arabic{Theorem}.\arabic{Claim}}
\section{Fairly large polarization of a finitely generated field}
\label{sec:suff:large:pol:fgf}
Let $K$ be a finitely generated field over $\QQ$ with
$d = \trdeg_{\QQ}(K)$, and let
$B$ be a flat and projective integral scheme over $\ZZ$ such that
$K$ is the function field of $B$. 
Let $\overline{L}$ be a $C^{\infty}$-hermitian $\QQ$-line bundle on $B$.
Here we fix several notations.

$\bullet${\bf nef}:
We say  $\overline{L}$ is {\em nef} if
$c_1(\overline{L})$ is a semipositive form on $B(\CC)$ and,
for all one-dimensional integral closed subschemes $\Gamma$ of $B$,
$\adeg \left( \rest{\overline{L}}{\Gamma} \right) \geq 0$.

$\bullet${\bf big}:
$\overline{L}$ is said to be {\em big} if 
$\rank_{\ZZ} H^0(B, L^{\otimes m}) = O(m^d)$, and there is a non-zero
section $s$ of $H^0(B, L^{\otimes n})$ with $\Vert s \Vert_{\sup} < 1$ for
some positive integer $n$.

$\bullet${\bf $\pmb{\QQ}$-effective}:
$\overline{L}$
is said to be {\em $\QQ$-effective} if there is a positive integer $n$ and
a non-zero $s \in H^0(B, L^{\otimes n})$ with $\Vert s \Vert_{\sup} \leq 1$.
If $\overline{L}_1 \otimes \overline{L}_2^{\otimes -1}$ is $\QQ$-effective
for $C^{\infty}$-hermitian $\QQ$-line bundles $\overline{L}_1, \overline{L}_2$
on $B$, then we denote this by $\overline{L}_1 \succsim \overline{L}_2$.

$\bullet${\bf polarization}:
A collection $\overline{B} = (B; \overline{H}_1, \ldots, \overline{H}_d)$
of $B$ and nef $C^{\infty}$-hermitian $\QQ$-line bundles 
$\overline{H}_1, \ldots,
\overline{H}_d$ on $B$ is called a {\em polarization of $K$}.

$\bullet${\bf fairly large polarization}:
A polarization $\overline{B} = (B; \overline{H}_1, \ldots, \overline{H}_d)$
is said to be {\em fairly large}
if there are generically finite morphisms
$\mu : B' \to B$ and $\nu : B' \to \left(\PP^1_{\ZZ}\right)^d$ of
flat and projective integral schemes over $\ZZ$, and
nef and big $C^{\infty}$-hermitian $\QQ$-line bundles
$\overline{L}_1, \ldots, \overline{L}_d$
on $\PP^1_{\ZZ}$ such that
$\mu^*(\overline{H}_i) \succsim 
\nu^*(p_i^*(\overline{L}_i))$ for all $i$,
where $p_i : \left(\PP^1_{\ZZ}\right)^d \to \PP^1_{\ZZ}$ is the projection
to the $i$-th factor.

\medskip
Finally we would like to give a simple sufficient condition for the fair
largeness of a polarization.
Let $k$ be a number field, and $O_k$ the ring of integer in $k$.
Let $B_1, \ldots, B_l$ be projective and flat integral schemes over $O_k$ whose
generic fibers over $O_k$ are geometrically irreducible.
Let $K_i$ be the function field of $B_i$ and $d_i$ 
the transcendence degree of $K_i$
over $k$. We set $B = B_1 \times_{O_k} \cdots \times_{O_k} B_l$ and
$d = d_1 + \cdots + d_l$. Then, the function field of $B$ is
the quotient field of $K_1 \otimes_k K_2 \otimes_k \cdots \otimes_k K_l$,
which is denoted by $K$, and
the transcendence degree of $K$ over $k$ is $d$.
For each $i$ ($i = 1, \ldots, l$),
let $\overline{H}_{i, 1}, \ldots, \overline{H}_{i, d_i}$
be nef and big $C^{\infty}$-hermitian $\QQ$-line bundles on $B_i$.
We denote by $q_i$ the projection
$B \to B_i$ to the $i$-th factor.
Then, we have the following.

\begin{Proposition}
\label{prop:large:polarization}
A polarization $\overline{B}$ of $K$ given by
\[
\overline{B} = \left(B; q_1^*(\overline{H}_{1,1}), \ldots, 
q_1^*(\overline{H}_{1, d_1}),
\ldots, q_l^*(\overline{H}_{l,1}), \ldots, q_l^*(\overline{H}_{l, d_l})
\right)
\]
is fairly large.
\end{Proposition}

\Proof
Since there is a dominant rational map 
$B_i \dashrightarrow \left( \PP^1_{\ZZ} \right)^{d_i}$
by virtue of Noether's normalization theorem,
we can find a birational morphism
$\mu_i : B'_i \to B_i$ of projective integral schemes over $O_k$ and
a generically finite morphism $\nu_i : B'_i \to  
\left( \PP^1_{\ZZ} \right)^{d_i}$.
We set $B' = B'_1 \times_{O_k} \cdots \times_{O_k} B'_l$,
$\mu = \mu_1 \times \cdots \times \mu_l$ and 
$\nu = \nu_1 \times \cdots \times \nu_l$.
Let $\overline{L}$ be a $C^{\infty}$-hermitian line bundle on 
$\PP^1_{\ZZ}$ given by
$(\OO_{\PP^1_{\ZZ}}(1), \Vert\cdot\Vert_{FS})$.
Note that $\overline{L}$ is nef and big.
Then, since $\mu_i^*(\overline{H}_{i,j})$ is big,
there is a positive integer $a_{i,j}$ with
$\mu_i^*(\overline{H}_{i,j})^{\otimes a_{i,j}} \succsim 
\nu_i^*\left(  p_j^*(\overline{L}) \right)$, that is,
$\mu_i^*(\overline{H}_{i,j}) \succsim 
\nu_i^*\left(  p_j^*\left(\overline{L}^{\otimes 1/a_{i,j}}\right) \right)$.
Thus, we get our proposition.
\QED

\section{Comparison of heights with respect to different polarizations}
First of all, let us recall the definition of height functions over
finitely generated fields.

Let $K$ be a finitely generated field over $\QQ$ with
$d = \trdeg_{\QQ}(K)$, and let
$\overline{B} = (B; \overline{H}_1, \ldots, \overline{H}_d)$
be a polarization of $K$.
First of all, let us recall the definition of height
of $\overline{K}$-valued points.
Let $X$ be a geometrically irreducible
projective variety over $K$ and $L$ an ample line bundle on $X$.
Let us take a projective integral scheme $\XX$ over $B$ and
a $C^{\infty}$-hermitian $\QQ$-line bundle $\overline{\LL}$ on $\XX$ such that
$X$ is the generic fiber of $\XX \to B$ and $L$ is equal to $\LL_K$
in $\Pic(X) \otimes \QQ$. Then, for $x \in X(\overline{K})$,
we define $h^{\overline{B}}_{(\XX, \LL)}(x)$ to be
\[
h^{\overline{B}}_{(\XX, \overline{\LL})}(x) =
\frac{\adeg \left( \acherncl_1(\overline{\LL}) \cdot
\prod_{j=1}^d \acherncl_1(\pi^*(\overline{H}_j)) \cdot \Delta_x \right)}
{[K(x) : K]},
\]
where $\Delta_x$ is the Zariski closure in $\XX$ of the image
$\Spec(\overline{K}) \to X \hookrightarrow \XX$, and
$\pi : \XX \to B$ is the canonical morphism.
By virtue of \cite{MoArht},
if $(\XX', \LL')$ is another model of $(X, L)$ over $B$, then
there is a constant $C$ with
$\vert h^{\overline{B}}_{(\XX, \LL)}(x) - 
h^{\overline{B}}_{(\XX', \LL')}(x) \vert \leq C$
for all $x \in X(\overline{K})$.
Hence, we have the unique height function $h^{\overline{B}}_L$
modulo the set of bounded functions.

More generally, we can define the height of cycles on $X_{\overline{K}}$.
We assume that $\overline{\LL}$ is nef with respect to $\pi : \XX
\to B$, that is,
\begin{enumerate}
\renewcommand{\labelenumi}{(\arabic{enumi})}
\item
For any analytic maps $h : M \to \XX(\CC)$ from a complex manifold $M$
to $\XX(\CC)$ with $\pi(h(M))$ being a point,
$c_1(h^*(\overline{\LL}))$ is semipositive.

\item
For every $b \in B$, the restriction $\rest{\LL}{X_{\overline{b}}}$ of $\LL$
to the geometric fiber over $b$ is nef.
\end{enumerate}
Let $Z$ be an effective cycle on $X_{\overline{K}}$.
We assume that $Z$ is defined over a finite extension field $K'$ of $K$.
Let $B'$ be
the normalization of $B$ in $K'$, and
let $\rho : B' \to B$ be the induced morphism.
Let $\XX'$ be the main component of $\XX \times_{B} B'$. 
We set the induced morphisms
as follows.
\[
\begin{CD}
\XX @<{\tau}<< \XX' \\
@V{\pi}VV @VV{\pi'}V \\
B     @<{\rho}<< B'
\end{CD}
\]
Let $\ZZZ$ be the Zariski closure of $Z$ in $\XX'$.
Then the height $h^{\overline{B}}_{(\XX, \overline{\LL})}(Z)$ of $Z$
with respect to $(\XX, \overline{\LL})$ and $\overline{B}$
is defined by
\[
h^{\overline{B}}_{(\XX, \overline{\LL})}(Z) 
= \frac{\adeg\left( 
\acherncl_1\left({\tau}^*(\overline{\LL})\right)^{\cdot \dim Z  +1} \cdot 
\prod_{j=1}^d
\acherncl_1\left({\pi'}^*({\rho}^*(\overline{H}_j))\right) \cdot \ZZZ
\right)}{[K':K](\dim Z + 1)\deg_L(Z)}.
\]
Note that the above definition does not depend on the choice of $K'$
by the projection formula.
Let $(\YY, \overline{\MM})$ be another model of $(X, L)$ over $B$
such that $\overline{\MM}$ is nef with respect to $\YY \to B$.
Then, there is a constant $C$ such that
\[
\left\vert 
h^{\overline{B}}_{(\XX, \overline{\LL})}(Z) - 
h^{\overline{B}}_{(\YY, \overline{\MM})}(Z)
\right\vert \leq C
\]
for all effective cycles $Z$ of $X_{\overline{K}}$
(cf. \cite[Proposition~2.1]{MoCanht}).
Thus, we may denote by $h^{\overline{B}}_L$
the class of $h^{\overline{B}}_{(\XX, \overline{\LL})}$ modulo
the set of bounded functions.
Moreover, we say $h^{\overline{B}}_L$ is the height function
associated with $L$ and $\overline{B}$.

\medskip
Let $k$ be a number field and $O_k$ the ring of integers in $k$.
Fix a positive integer $d$. For each $1 \leq i \leq d$, let
$B_i$ be a flat and projective integral scheme over $O_k$ whose generic fiber
over $O_k$ is a geometrically irreducible curve over $k$.
Let $\overline{M}_i$ be a nef and big hermitian $\QQ$-line bundle on $B_i$.
Moreover, let $B$ be a flat and projective integral scheme over $O_k$, and
$\nu : B \to B_1 \times_{O_k} \cdots \times_{O_k} B_d$ a
generically finite morphism.
We denote the function field of $B$ (resp. $B_i$) by $K$ (resp. $K_i$).
Note that $\trdeg_{k}(K) = d$ and $\trdeg_{k}(K_i) = 1$ for all $i$.
We set $\overline{H}_i = \nu^*p_i^*(\overline{M}_i)$ for each $i$ and
$\overline{H} = \bigotimes_{i=1}^d \overline{H}_i$,
where $p_i : B_1 \times_{O_k} \cdots \times_{O_k} B_d \to B_i$
is the projection to the $i$-th factor.
Further, we set
\[
\lambda_i = \exp\left(-\frac{\adeg(\acherncl_1(\overline{M}_i)^2)}{
[k : \QQ]\deg((M_i)_k)}\right).
\]
Here we consider several kinds of polarizations of $K$ as follows:
\[
\begin{cases}
\overline{B}_0 = (B; \overline{H}, \ldots, \overline{H}), \\
\overline{B}_1 = (B; \overline{H}_1, \ldots, \overline{H}_d), \\
\overline{B}_{i,j} = (B; \overline{H}_1, \ldots, \overline{H}_{j-1},
(\OO_B, \lambda_i\vert\cdot\vert_{can}), \overline{H}_{j+1}, \ldots,
\overline{H}_d) & \text{for $i \not= j$}.
\end{cases}
\]
A key result of this note is the following.

\begin{Proposition}
\label{prop:comp:heights}
Let $X$ be a geometrically irreducible projective variety over $K$, and
$L$ an ample line bundle on $X$. Let $(\XX, \overline{\LL})$ be
a model of $(X, L)$ over $B$ such that
$\overline{\LL}$ is nef with respect to $\XX \to B$. 
Then, for all effective cycles $Z$ on $X_{\overline{K}}$,
\[
h^{\overline{B}_0}_{(\XX, \overline{\LL})}(Z)
= d! h^{\overline{B}_1}_{(\XX, \overline{\LL})}(Z) + \frac{d!}{2}
\sum_{i\not=j} h^{\overline{B}_{i,j}}_{(\XX, \overline{\LL})}(Z).
\]
In particular, we can find a constant $C$ such that
\[
h^{\overline{B}_0}_{L}(Z)
\leq C h^{\overline{B}_1}_{L}(Z) + O(1)
\]
for all effective cycles $Z$ on $X_{\overline{K}}$.
\end{Proposition}

\Proof
Let $Z$ be an effective cycle on $X_{\overline{K}}$.
We assume that $Z$ is defined over a finite extension field $K'$ of $K$.
Let $B'$ be
the normalization of $B$ in $K'$, and
let $\rho : B' \to B$ be the induced morphism.
Let $\XX'$ be the main component of $\XX \times_{B} B'$. 
We set the induced morphisms
as follows.
\[
\begin{CD}
\XX @<{\tau}<< \XX' \\
@V{\pi}VV @VV{\pi'}V \\
B     @<{\rho}<< B'
\end{CD}
\]
Let $\ZZZ$ be the Zariski closure of $Z$ in $\XX'$.
We set $\overline{D} =
\nu_* \rho_*\pi'_*\left(\acherncl_1\left(
{\tau}^*(\overline{\LL})\right)^{\cdot \dim Z  +1}
\cdot \ZZZ\right)$.
Then,
\begin{align*}
h^{\overline{B}_0}_{(\XX, \overline{\LL})}(Z) & = 
\frac{\adeg\left(\overline{D} \cdot \left( \sum_{l=1}^d \acherncl_1
\left(p_l^*(\overline{M}_l)\right) \right)^{\cdot d}
 \right)}{[K':K](\dim Z + 1)\deg_L(Z)} \\
& = \sum_{a_1 + \cdots + a_d = d}
\frac{d!}{a_1! \cdots a_d!}
\frac{\adeg\left(\overline{D} \cdot \prod_{l=1}^d 
\acherncl_1(p_l^*(\overline{M}_l))^{\cdot a_l}
\right)}{[K':K](\dim Z + 1)\deg_L(Z)}.
\end{align*}
Here we claim the following.

\begin{Claim}
If $(a_1, \ldots, a_d) \not= (1, \ldots, 1)$ and
$\adeg\left(\overline{D} \cdot \prod_{l=1}^d 
\acherncl_1(p_l^*(\overline{M}_l))^{\cdot a_l}
\right) \not= 0$, then there are $i, j \in \{1, \ldots, d\}$ such that
$a_i = 2$, $a_j= 0$ and $a_l = 1$ for all $l \not= i, j$. In particular,
\[
h^{\overline{B}_0}_{(\XX, \overline{\LL})}(Z) =
d! h^{\overline{B}_1}_{(\XX, \overline{\LL})}(Z) + \frac{d!}{2}
\sum_{i\not=j} \frac{\adeg\left(\overline{D} \cdot 
\acherncl_1(p_i^*(\overline{M}_i))^{\cdot 2}
\prod_{l=1, l\not=i,j}^d \acherncl_1(p_l^*(\overline{M}_l))
\right)}{[K':K](\dim Z + 1)\deg_L(Z)}.
\]
\end{Claim}

Clearly, $a_l \leq 2$ for all $l$. Thus, there is $i$ with $a_i = 2$.
Suppose that $a_j = 2$ for some $j \not= i$.
Then,
\[
\adeg\left(\overline{D} \cdot \prod_{l=1}^d 
\acherncl_1(p_l^*(\overline{M}_l))^{\cdot a_l}
\right) = \adeg\left(\overline{D} \cdot 
\acherncl_1(p_i^*(\overline{M}_i))^{\cdot 2} \cdot
\acherncl_1(p_j^*(\overline{M}_j))^{\cdot 2} \cdot
\prod_{l=1, l\not=i,j}^d \acherncl_1(p_l^*(\overline{M}_l))^{\cdot a_l}
\right).
\]
Thus, using the projection formula with respect to $p_i$,
\[
\adeg\left(\overline{D} \cdot \prod_{l=1}^d 
\acherncl_1(p_l^*(\overline{M}_l))^{\cdot a_l}
\right) =
\adeg(\acherncl_1(\overline{M}_i)^{\cdot 2})
\deg \left( D_{\eta_i} \cdot
p_j^*(M_j)_{\eta_i}^{\cdot 2} \cdot 
\prod_{l=1, l\not=i,j}^d p_l^*(\overline{M}_l)^{\cdot a_l}_{\eta_i} \right),
\]
where $\eta_i$ means the restriction to the generic fiber of $p_i$.
Here the generic fiber of $p_i$ is isomorphic to
$(B_1 \times_{k} K_i) \times_{K_i} \cdots
(B_{i-1} \times_{k} K_i) \times_{K_i} (B_{i+1} \times_{k} K_i) \times
(B_{d} \times_{k} K_i)$ and $B_j \times_k K_i$ is a projective curve over $K_i$.
Thus, we can see
\[
\adeg\left(\overline{D} \cdot \prod_{l=1}^d 
\acherncl_1(p_l^*(\overline{M}_l))^{\cdot a_l}
\right) = 0.
\]
This is a contradiction.
Hence, we get our claim.

\medskip
By the above claim, it is sufficient to see that
\addtocounter{Claim}{1}
\begin{equation}
\label{eqn:prop:comp:heights:10}
h^{\overline{B}_{i,j}}_{(\XX, \overline{\LL})}(Z) =
\frac{\adeg\left(\overline{D} \cdot 
\acherncl_1(p_i^*(\overline{M}_i))^{\cdot 2}
\prod_{l=1, l\not=i,j}^d \acherncl_1(p_l^*(\overline{M}_l))
\right)}{[K':K](\dim Z + 1)\deg_L(Z)}.
\end{equation}
First of all,
\[
h^{\overline{B}_{i,j}}_{(\XX, \overline{\LL})}(Z) =
\frac{{\displaystyle -\log(\lambda_i) \int_{\ZZZ(\CC)} 
c_1(\tau^* \overline{\LL})^{\wedge \dim Z +1} \wedge
\bigwedge_{l=1,l\not=j} c_1({\pi'}^* \rho^* \nu^*p_l^*(\overline{M}_l))}}
{[K':K](\dim Z + 1)\deg_L(Z)}.
\]
Moreover,
\[
\int_{\ZZZ(\CC)} c_1(\tau^*\overline{\LL})^{\wedge \dim Z +1} \wedge
\bigwedge_{l=1,l\not=j} c_1({\pi'}^* \rho^*\nu^*p_l^*(\overline{M}_l))
\]
is equal to
\[
[k : \QQ] \deg\left(\tau^* \LL_k^{\cdot \dim Z + 1} \cdot
\prod_{l=1,l\not=j} {\pi'}^* \rho^*\nu^*p_l^*(\overline{M}_l)_k
\cdot \ZZZ_k \right) =
[k : \QQ] \deg\left(D_k \cdot \prod_{l=1,l\not=j} p_l^*(\overline{M}_l)_k
\right).
\]
Thus, we obtain
\[
h^{\overline{B}_{i,j}}_{(\XX, \overline{\LL})}(Z) 
= \frac{\adeg(\acherncl_1(\overline{M}_i)^2)\deg\left(D_k \cdot
\prod_{l=1,l\not=j} p_l^*(\overline{M}_l)_k\right)}{
\deg((M_i)_k)[K':K](\dim Z + 1)\deg_L(Z)}.
\]
On the other hand, by the projection formula with respect to $p_i$,
\[
\adeg\left(\overline{D} \cdot 
\acherncl_1(p_i^*(\overline{M}_i))^{\cdot 2}
\prod_{l=1, l\not=i,j}^d \acherncl_1(p_l^*(\overline{M}_l))
\right) = \adeg(\acherncl_1(\overline{M}_i)^{\cdot 2}) \cdot
\deg \left( D_{\eta_i} \prod_{l=1, l\not=i,j}^d 
p_l^*(\overline{M}_l)_{\eta_i} \right),
\]
where $\eta_i$ means the restriction to the generic fiber of $p_i$.
Moreover, by the projection formula again,
\[
\deg\left(D_k \cdot \prod_{l=1,l\not=j} p_l^*(\overline{M}_l)_k
\right) = \deg(M_i)_k \cdot \deg \left( D_{\eta_i} 
\prod_{l=1, l\not=i,j}^d p_l^*(\overline{M}_l)_{\eta_i}
\right).
\]
Thus, we get \eqref{eqn:prop:comp:heights:10}.

The last assertion is obvious because
there is a positive integer $m$ such that
\[
\overline{H}_j^{\otimes m} \succsim (\OO_B, \lambda_i\vert\cdot\vert_{can})
\]
for every $i,j$.
\QED

\begin{Corollary}
\label{cor:comp:big:large:pol}
Let $X$ be a geometrically irreducible projective variety over $K$, and
$L$ an ample line bundle on $X$. Let
$\overline{B}$ and $\overline{B}'$ be polarizations of $K$.
We assume that $\overline{B}$ is big and $\overline{B}'$ is fairly large.
Then, there are positive constants $a$ and $b$ such that
\[
a h_L^{\overline{B}'}(Z) + O(1) \leq h_L^{\overline{B}}(x)
\leq b h_L^{\overline{B}'}(Z) + O(1)
\]
for all effective cycles $Z$ on $X_{\overline{K}}$.
In particular, if $X$ is an abelian variety and
$L$ is a symmetric ample line bundle,
then 
\[
a \hat{h}_L^{\overline{B}'}(Z) \leq \hat{h}_L^{\overline{B}}(Z)
\leq b \hat{h}_L^{\overline{B}'}(Z)
\]
for all effective cycles $Z$ on $X_{\overline{K}}$.
\end{Corollary}

\Proof
The first inequality is a consequence of 
\cite[(5) of Proposition~3.3.7]{MoArht}.
We set $\overline{B} = (B; \overline{H}_1, \ldots, \overline{H}_d)$ and
$\overline{B}' = (B'; \overline{H}'_1, \ldots, \overline{H}'_d)$.
Since $\overline{B}'$ is fairly large,
there are generically finite morphisms
$\mu' : B'' \to B'$ and $\nu : B'' \to \left(\PP^1_{\ZZ}\right)^d$ of
flat and projective integral schemes over $\ZZ$,
and nef and big $C^{\infty}$-hermitian $\QQ$-line bundles
$\overline{L}_1, \ldots, \overline{L}_d$
on $\PP^1_{\ZZ}$ such that
${\mu'}^*(\overline{H}'_i) \succsim 
\nu^*(p_i^*(\overline{L}_i))$ for all $i$,
where $p_i : \left(\PP^1_{\ZZ}\right)^d \to \PP^1_{\ZZ}$ is the projection
to the $i$-th factor. Changing $B''$ if necessarily, we may assume that there is
a generically finite morphism $\mu : B'' \to B$.
By virtue of the projection formula, we may assume that
$B = B' = B''$.
We set $\overline{H} = \nu^*\left( \bigotimes_{l=1}^d p_l^*(\overline{L}_l) 
\right)$. 
Then, $(B; \overline{H},\ldots, \overline{H})$ is a big polarization. 
Thus, there is a positive integer $b_1$ such that
\[
h_L^{\overline{B}} \leq b_1 h_L^{(B; \overline{H}, \ldots, \overline{H})} 
+ O(1).
\]
Moreover, by Proposition~\ref{prop:comp:heights}, 
we can find a positive constant $b_2$
with
\[
h_L^{(B; \overline{H}, \ldots, \overline{H})} \leq b_2 
h_L^{(B; \nu^*p_1^*(\overline{L}_1), \ldots, 
\nu^*p_d^*(\overline{L}_d))} + O(1).
\]
On the other hand, since ${\overline{H}'_i} \succsim 
\nu^*(p_i^*(\overline{L}_i))$ for all $i$,
\[
h_L^{(B; \nu^*p_1^*(\overline{L}_1), \ldots, \nu^*p_d^*(\overline{L}_d))}
\leq h_L^{(B; \overline{H}_1, \ldots, \overline{H}_d)} + O(1).
\]
Hence, we get our corollary.
\QED

Here, let us give the proof of Theorem~\ref{thm:northcott:intro},
Theorem~\ref{thm:BL:conj:fair:pol:intro} and
Theorem~\ref{thm:bogo:equiv:height} in the
introduction. 
Theorem~\ref{thm:northcott:intro} is obvious by
\cite[Theorem~4.3]{MoArht} and
Corollary~\ref{cor:comp:big:large:pol}. 
Theorem~\ref{thm:BL:conj:fair:pol:intro} is a consequence of
\cite{MoBL}, Corollary~\ref{cor:comp:big:large:pol} and
the following lemma.

\begin{Lemma}
Let $V$ be a vector space over $\RR$, and
$\langle\ , \ \rangle$ and $\langle\ , \ \rangle'$ be two inner products on $V$.
If $\langle x , x \rangle \leq \langle x , x \rangle'$ for all $x \in V$, then
$\det \left( \langle x_i , x_j \rangle \right) \leq \det 
\left( \langle x_i , x_j \rangle' \right)$
for all $x_1, \ldots, x_n \in V$.
\end{Lemma}

\Proof
If $x_1, \ldots, x_n$ are linearly dependent, then our assertion is trivial.
Otherwise, it is nothing more than \cite[Lemma~3.4]{MoBG}.
\QED

Finally, let us consider the proof of Theorem~\ref{thm:bogo:equiv:height}.
The equivalence of (1) and (2) follows from
Theorem~\ref{thm:BL:conj:fair:pol:intro}.
It is obvious that (1) implies (3).
Conversely, (3) implies (1) by virtue of
\cite{MoCanht} and Corollary~\ref{cor:comp:big:large:pol}.

\section{Northcott's theorem for cycles}
\label{sec:def:height:cycle}
In this section, 
we will generalize Northcott's theorem to 
the height of cycles on projective varieties.

\begin{Theorem}
\label{thm:northcott:cycle:intro}
Let $K$ be a finitely generated field over $\QQ$,
$\overline{K}$ the algebraic closure of $K$,
and let $\overline{B}$ be a fairly large
polarization of $K$.
Let $X$ be a geometrically irreducible projective variety over $K$, and
$L$ an ample line bundle on $X$.
For an effective cycle $Z$ on $X_{\overline{K}}$,
we denote by $h^{\overline{B}}_L(Z)$ the height of $Z$ 
with respect to $\overline{B}$.
Moreover, the orbit of $Z$ by the action of the Galois group 
$\Gal(\overline{K}/K)$ is denoted
by $O_{\Gal(\overline{K}/K)}(Z)$.
Then, for a real number $M$ and integers $l$ and $e$, the set
of all effective cycles on 
$X_{\overline{K}}$ with
$h^{\overline{B}}_L(Z) \leq M$, $\deg_L(Z) \leq l$ and
$\#O_{\Gal(\overline{K}/K)}(Z)\leq e$
is finite.
\end{Theorem}

\Proof
Let us begin with the following lemma.

\begin{Lemma}
\label{lem:comp:norms}
Let $X = \PP_{\CC}^{n_1} \times \cdots \times \PP_{\CC}^{n_r}$ be
a product of projective spaces and
$\OO(d_1, \ldots, d_r)$ the line bundle on $X$
with the multi-degree $(d_1, \ldots, d_r)$.
Let $\Vert\cdot\Vert$ be the hermitian metric of
$\OO(d_1, \ldots, d_r)$ given by the Fubini-Study metric.
For each $i$, we fix a basis of
$H^0(\PP^{n_i}_{\CC}, \OO_{\PP^{n_i}_{\CC}}(1))$.
Then $H^0(X, \OO(d_1, \ldots, d_r))$ is
naturally isomorphic to the space of homogeneous
polynomials with the multi-degree $(d_1, \ldots, d_r)$.
For $s \in  H^0(X, \OO(d_1, \ldots, d_r))$,
we denote by $\vert s \vert$ the maximal value of the absolute of
coefficient of $s$ as a polynomial.
Then, there is a constant $C$ depending only on
$n_1, \ldots, n_r$, $d_1, \ldots, d_r$ and
a basis of $H^0(\PP^{n_i}_{\CC}, \OO_{\PP^{n_i}_{\CC}}(1))$
for each $i$ such that
\[ 
\vert s \vert \leq C \exp \left(
\int_X \Vert s \Vert \bigwedge_{j=1}^{r} 
p_j^*( c_1(\OO(1), \Vert\cdot\Vert_{FS}))^{\wedge n_j}
\right)
\]
for all $s \in  H^0(X, \OO(d_1, \ldots, d_r))$.
\end{Lemma}

\Proof
By virtue of \cite[Corollary~1.4.3]{BGS},
\[
\sup_{x \in X} \{ \Vert s \Vert(x) \}
\leq
 \exp\left( \sum_{i=1}^r \sum_{m=1}^{n_i} 
\frac{d_i}{2m} \right) \cdot \exp\left(
\int_X \Vert s \Vert \bigwedge_{j=1}^r 
p_j^*( c_1(\OO(1), \Vert\cdot\Vert_{FS}))^{\wedge n_j}
\right).
\]
On the other hand, $\sup_{x \in X} \{ \Vert s \Vert(x) \}$ and
$\vert s \vert$ give rise to two norms on the
finite dimensional space $ H^0(X, \OO(d_1, \ldots, d_r))$.
Thus, we get our lemma.
\QED

Let us start the proof of Northcott's theorem for cycles.
It is sufficient to see that, for any real number $M$ and 
any integers $e$, $d$ and $l$,
the set
\[
\left\{ Z \left|
\begin{array}{l}
\text{$Z$ is an effective cycle with} \\
\text{$\#O_{\Gal(\overline{K}/K)}(Z) \leq e$, $\deg_L(Z) = d$,
$\dim Z = l$ and $h^{\overline{B}}_L(Z) \leq M$} 
\end{array}
\right.\right\}
\]
is finite. Clearly, we may assume that $X = \PP^n_K$ and $L = \OO_{\PP^n_K}(1)$.
Let $K'$ be the invariant field of the stabilizer at $Z$. Then,
$[K' : K] \leq e$ and $Z$ is defined over $K'$.
Let $B'$ be the normalization of $B$ in $K'$, and
let $\overline{H}'_1, \ldots, \overline{H}'_d$ be the pull-backs of
$\overline{H}_1, \ldots, \overline{H}_d$ by
$B' \to B$ respectively.
Let $\ZZZ$ be the Zariski closure of $Z$ in 
$\PP^n_{B'} = \PP^n_{\ZZ} \times B'$.
Then,
\[
h^{\overline{B}}_L(Z) =
\frac{\adeg\left(
\acherncl_1(\overline{\OO_{\PP^n_{B'}}(1)})^{l+1} \cdot \prod_{j=1}^d
\acherncl_1(\pi^*_{B'}(\overline{H}'_j)) \cdot \ZZZ \right)}
{[K' : K] (l+1) \deg_L(Z)} +O(1),
\]
where $\pi_{B'}$ is the canonical projection $\PP^n_{B'} \to B'$ and
$\overline{\OO_{\PP^n_{B'}}(1)}$ is the pull-back of
$(\OO_{\PP^n_{\ZZ}}(1), \Vert\cdot\Vert_{FS})$ via $\PP^n_{B'} \to \PP^n_{\ZZ}$.

Let $\Check{\PP}^n_{\ZZ}$ be the dual projective space of $\PP^n_{\ZZ}$.
Let us consider
\[
\left( \Check{\PP}^n_{B'} \right)^{\times_{B'} (l+1)} =
\overbrace{\Check{\PP}_{B'}^n \times_{B'} \cdots \times_{B'} 
\Check{\PP}_{B'}^n}^{l+1},
\]
where $\Check{\PP}^n_{B'} = \Check{\PP}^n_{\ZZ} \times B'$.
Let $p_i : \left( \Check{\PP}^n_{B'} \right)^{\times_{B'} (l+1)}  \to
\Check{\PP}^n_{B'}$ be the projection to the $i$-th factor and
$p : \left( \Check{\PP}^n_{B'} \right)^{\times_{B'} (l+1)} \to B'$
the canonical morphism.
We set
\[
\OO_{B'}(d, \ldots, d) = \bigotimes_{i=1}^{l+1} 
p_i^*(\OO_{ \Check{\PP}^n_{B'}}(1)).
\]
Let $\Chow(Z)$ be the Chow divisor of $Z$, i.e.,
$\Chow(Z)$ is an element of $\vert \OO_{K'}(d, \ldots, d) \vert$
on $\left(\Check{\PP}_{K'}^n\right)^{l+1}$.
Let $\Chow(\ZZZ)$ be the Zariski closure of
$\Chow(Z)$ in $\left(\Check{\PP}_{B'}^n\right)^{\times_{B'}(l+1)}$.
Here we claim the following equation:
\addtocounter{Claim}{1}
\begin{multline}
\label{eqn:thm:northcott:cycle:1}
p_*\left( \prod_{i=1}^{l+1} \acherncl_1
\left(p_i^*(\overline{\OO_{\Check{\PP}^n_{B'}}(1)})\right)^{\cdot n} 
\cdot \Chow(\ZZZ) \right) \\
=
\pi_*\left(\acherncl_1(\overline{\OO_{\PP^n_{B'}}(1)})^{\cdot l+1} \cdot \ZZZ
\right) +
d(l+1) \pi_*
\left(\acherncl_1(\overline{\OO_{\PP^n_{B'}}(1)})^{\cdot n+1}\right).
\end{multline}
Let $U$ be the maximal Zariski open set of $B'$ such that
$\ZZZ \to B'$ is flat over $U$.
Then, in the same way as in the proof of
\cite[Proposition~1.2]{BSS}, we can see that
the equation \eqref{eqn:thm:northcott:cycle:1} holds over $U$.
Therefore, so does over $B'$ by \cite[Lemma~2.5.1]{KMSemi} because
$\codim(B' \setminus U) \geq 2$.

\medskip
By using \eqref{eqn:thm:northcott:cycle:1}, we can see
\addtocounter{Claim}{1}
\begin{multline}
\label{eqn:thm:northcott:cycle:2}
\left(h^{\overline{B}}_L(Z) +O(1)\right)
+ \adeg(\acherncl_1(\overline{\OO_{\PP^n_{\ZZ}}(1)})) \deg({H_1}_{\QQ} \cdots
{H_d}_{\QQ}) \\
= \frac{
\adeg \left(\prod_{i=1}^{l+1}\acherncl_1
\left(p_i^*(\overline{\OO_{\Check{\PP}^n_{B'}}(1)})\right)^{\cdot n} 
 \cdot \prod_{j=1}^d \acherncl_1(p^*(\overline{H}'_j)) \cdot
\Chow(\ZZZ) \right)}{[K' : K]}. 
\end{multline}
Choose $P \in H^0\left(\left(\Check{\PP}_{K'}^n\right)^{l+1}, 
\OO_{K'}(d, \ldots, d) \right)$
with $\zero(P) = \Chow(Z)$. 
Here, we fix a basis of $H^0(\PP^n_{\ZZ}, \OO_{\PP^n_{\ZZ}}(1))$.
Then, $P$ can be written by a polynomial
with coefficients $\{ a_{\lambda} \}_{\lambda \in \Lambda}$ in $K'$, that is,
$\{ a_{\lambda} \}_{\lambda \in \Lambda}$ is a Chow coordinate of $Z$.
Noting that $P$ gives rise a rational section of 
$\OO_{B'}(d,\ldots, d)$ over $B'$,
let
\[
\zero(P) = \Chow(\ZZZ) + \sum_{\Gamma} c_{\Gamma} p^{*}(\Gamma)
\]
be the decomposition as a rational section of $\OO_{B'}(d,\ldots, d)$,
where $\Gamma$ runs over all prime divisors on $B'$.
Then, we can easily see $c_{\Gamma} = 
\min_{\lambda \in \Lambda} \{ \ord_{\Gamma}(a_{\lambda}) \}$.
Here, let us calculate 
\[
\adeg \left(\prod_{i=1}^{l+1} \acherncl_1\left(
p_i^*(\overline{\OO_{\Check{\PP}^n_{B'}}(1)})\right)^{\cdot n} 
\cdot \prod_{j=1}^d \acherncl_1(p^* \overline{H}'_j) \cdot
\acherncl_1(\overline{\OO_{B'}(d, \ldots, d)}) \right)
\]
in terms of the  rational section $P$.
Then,
\addtocounter{Claim}{1}
\begin{multline}
\label{eqn:thm:northcott:cycle:3}
\adeg \left( \prod_{i=1}^{l+1} \acherncl_1\left(
p_i^*(\overline{\OO_{\Check{\PP}^n_{B'}}(1)})\right)^{\cdot n}
\cdot \prod_{j=1}^d \acherncl_1(p^*(\overline{H}'_j)) \cdot
\Chow(\ZZZ) \right) = \\
[K' : K]
\adeg \left(\prod_{i=1}^{l+1} \acherncl_1\left(
p_i^*(\overline{\OO_{\Check{\PP}^n_{B}}(1)})\right)^{\cdot n}
\cdot \prod_{j=1}^d \acherncl_1(p^*(\overline{H}_j)) \cdot
\acherncl_1(\overline{\OO_B(d, \ldots, d)}) \right) \\
+ \sum_{\Gamma} \max_{\lambda \in \Lambda} \{ -\ord_{\Gamma}(a_{\lambda}) \}
\adeg \left( \prod_{j=1}^d \acherncl_1(\overline{H}'_j) \cdot
\Gamma \right)  \\
+ \int_{\left(\Check{\PP}_{\CC}^n\right)^{l+1} \times B'(\CC)}
\log \Vert P \Vert \bigwedge_{i=1}^{l+1}
c_1\left(p_j^*(\overline{\OO_{\Check{\PP}^n_{B'}}(1)})\right)^{\wedge n} \wedge
\bigwedge_{j=1}^d c_1(p^*(\overline{H}'_j)).
\end{multline}
Let $U_0$ be a Zariski open set of $B'(\CC)$ such that
$a_{\lambda}$ has no zeros and poles on $U_0$ for every $\lambda \in \Lambda$.
For each $b \in U_0$, let $i_b :  \left(\Check{\PP}_{\CC}^n\right)^{l+1} \to
\left(\Check{\PP}_{\CC}^n\right)^{l+1} \times B'(\CC)$ be a morphism
given by $i_b(x) = (x, b)$.
Then,
\begin{multline*}
\left(
\int_{\left(\Check{\PP}_{\CC}^n\right)^{l+1} \times B'(\CC)}
\log \Vert P \Vert \bigwedge_{i=1}^{l+1}
c_1\left(p_j^*(\overline{\OO_{\Check{\PP}^n_{B'}}(1)})\right)^{\wedge n}
\right)(b) \\
= \int_{\left(\Check{\PP}_{\CC}^n\right)^{l+1}}
i_b^*\left(\log \Vert P \Vert\right) \bigwedge_{i=1}^{l+1}
c_1\left(p_j^*(\overline{\OO_{\Check{\PP}^n_{\CC}}(1)})\right)^{\wedge n}
\end{multline*}
Thus, by Lemma~\ref{lem:comp:norms},
\addtocounter{Claim}{1}
\begin{multline}
\label{eqn:thm:northcott:cycle:4}
\int_{\left(\Check{\PP}_{\CC}^n\right)^{l+1} \times B'(\CC)}
\log \Vert P \Vert \bigwedge_{i=1}^{l+1}
c_1\left(p_j^*(\overline{\OO_{\Check{\PP}^n_{B'}}(1)})\right)^{\wedge n} \wedge
\bigwedge_{j=1}^d c_1(p^*(\overline{H}'_j)) \\
\geq \int_{B'(\CC)}
\log \max_{\lambda \in \Lambda} \{ \vert a_{\lambda} \vert \}
\bigwedge_{j=1}^d c_1(\overline{H}'_j) - [K':K] C'\int_{B(\CC)}
\bigwedge_{j=1}^d c_1(\overline{H}_j)
\end{multline}
for some constant $C'$ depending only on $n$, $d$, $l$ and
a basis $H^0(\PP^n_{\ZZ}, \OO_{\PP^n_{\ZZ}}(1))$.
Thus, gathering \eqref{eqn:thm:northcott:cycle:2}, 
\eqref{eqn:thm:northcott:cycle:3} and
\eqref{eqn:thm:northcott:cycle:4},
\[
h^{\overline{B}}_{nv}\left((a_{\lambda})_{\lambda \in \Lambda}\right) \leq
h^{\overline{B}}_L(Z) + C''
\]
for some constant $C''$ independent on $Z$.
Thus, we have only finitely many $(a_{\lambda})_{\lambda \in \Lambda}$ modulo
the scalar product of $\overline{K}$.
Therefore, we have only finitely many $\Chow(Z)$.
Hence we obtain our assertion because the correspondence
\[
\Chow :
\left\{ 
\begin{array}{l}
\text{Effective cycles $Z$ on $X_{\overline{K}}$ with} \\
\text{$l = \dim Z$ and $\deg(Z) = d$} 
\end{array}
\right\} \to
|\OO_{\overline{K}}(d, \ldots, d)|
\]
is injective.
\QED

\renewcommand{\theTheorem}{\Alph{section}.\arabic{Theorem}}
\renewcommand{\theClaim}{\Alph{section}.\arabic{Theorem}}
\renewcommand{\theequation}{\Alph{section}.\arabic{Theorem}}
\renewcommand{\thesection}{Appendix}
\setcounter{section}{0}
\setcounter{Theorem}{0}
\section{Direct proof of Northcott's theorem with respect to 
a fairly large polarization}
If we use a fairly large polarization, 
we can give a simpler proof of Northcott's theorem.
In this appendix, let us consider this problem.

Let us start the direct proof of Theorem~\ref{thm:northcott:intro}.
We denote by
$\overline{\OO_{\PP^1_{\ZZ}}(1)}$ the hermitian line bundle
$(\OO_{\PP^1_{\ZZ}}(1), \Vert\cdot\Vert_{FS})$ on $\PP^1_{\ZZ}$. 
First of all, we claim the following.

\addtocounter{Theorem}{1}
\begin{Claim}
We assume that there is a generically finite morphism
$\nu : B \to \left( \PP^1_{\ZZ} \right)^d$ with
$\overline{H}_i = \nu^*(p_i^*(\overline{\OO_{\PP^1_{\ZZ}}(1)}))$ for all $i$.
Then, the assertion of our theorem holds.
\end{Claim}

Since $L$ is ample, there is a positive integer $m$ and
an embedding $\phi : X \hookrightarrow \PP^n$ with 
$\phi^*(\OO_{\PP^n}(1)) = L^{\otimes m}$.
Thus, we may assume that $X = \PP^n_K$ and $L = \OO_{\PP^n_K}(1)$.
Let $\overline{B}_0$ be a polarization given by
$\left( \left( \PP^1_{\ZZ} \right)^d; p_1^*(\overline{\OO_{\PP^1_{\ZZ}}(1)}),
\ldots, p_d^*(\overline{\OO_{\PP^1_{\ZZ}}(1)}) \right)$.
Then, by the projection formula, we can see that
$h_L^{\overline{B}_0} = \deg(\nu) h_L^{\overline{B}} + O(1)$.
Therefore, we may assume $\overline{B} = \overline{B}_0$.
Moreover, in the same argument as in \cite[Claim~4.3.3]{MoArht}, 
we may assume $e=1$.

Let $\Delta_{\infty}$ be the closure of $\infty \in \PP^1_{\QQ}$
in $\PP^1_{\ZZ}$. We set $\Delta^{(i)}_{\infty} = p_i^*(\Delta_{\infty})$.
Moreover, we set $\overline{A}_i = 
p_i^*(\OO_{\PP^1_{\ZZ}}, (1/2)\vert\cdot\vert_{can})$.
If we denote $c_1(\PP^1, \Vert\cdot\Vert_{FS})$ by $\omega$, then
\addtocounter{Theorem}{1}
\begin{equation}
\label{eqn:thm:northcott:1}
\adeg\left(
\acherncl_1(\overline{H}_1) \cdots
\acherncl_1(\overline{A}_i) \cdots
\acherncl_1(\overline{H}_d) \cdot \Delta^{(j)}_{\infty}
\right) = \int_{p_j^{-1}(\infty)} \log(2)
\bigwedge_{l=1,l\not=i}^d p_l^*(\omega) = 
\begin{cases}
\log(2) & \text{if $j = i$} \\
0 & \text{if $j \not= i$}
\end{cases}
\end{equation}

Let $\overline{B}_i$ be a polarization of
$K$ given by $\overline{B}_i = (B;\overline{H}_1, \ldots,
\overline{H}_{i-1}, \overline{A}_i, \overline{H}_{i+1}, \ldots,
\overline{H}_d)$.
Here $\overline{\OO_{\PP^1_{\ZZ}}(1)}^{\otimes 2} \succsim
(\OO_{\PP^1_{\ZZ}}, (1/2)\vert\cdot\vert_{can})$ because
$\sup_{x \in \PP^1(\CC)} \Vert X_0 X_1 \Vert_{FS}(x) \leq 1/2$,
where $\{ X_0, X_1\}$ are a basis of
$H^0(\PP^1_{\ZZ},\OO_{\PP^1_{\ZZ}}(1))$.
Thus, by \cite[(5) of Proposition~3.3.7]{MoArht},
there are positive constants $a$ such that
$h_{nv}^{\overline{B}_i} \leq 2 h_{nv}^{\overline{B}} + a$ for all $i$.
We set
\[
 S =
\{ P \in \PP^n(\QQ(z_1, \cdots, z_d)) \mid h_{nv}^{\overline{B}}(P) \leq M \}
\]
Then, for any $P \in S$,
$h_{nv}^{\overline{B}_i}(P) \leq 2 M + a$.
Moreover, there are $f_0, \cdots, f_n \in \ZZ[z_1, \cdots, z_d]$
such that $f_0, \cdots, f_n$ are relatively prime and
$P = (f_0 : \cdots : f_n)$.
Here, by using \eqref{eqn:thm:northcott:1} together with
facts that $c_1(\overline{A}_i) = 0$ and
$f_0, \cdots, f_n$ are relatively prime,
\[
 h_{nv}^{\overline{B}_i}(P) =
\max \{ \deg_i(f_0), \ldots, \deg_i(f_n) \} \log(2),
\]
where $\deg_i$ is the degree of polynomials with respect to $z_i$.
Thus, there is a constant $M_1$ independent on $P \in S$
such that $\deg_i(f_j) \leq M_1$ for
all $i,j$.
On the other hand,
\begin{multline*}
 h_{nv}^{\overline{B}}(P) =
\sum_i \max \{ \deg_i(f_0), \ldots, \deg_i(f_n) \} 
\adeg\left( \acherncl_1(\overline{H}_1) \cdots
\acherncl_1(\overline{H}_d) \cdot \Delta^{(i)}_{\infty}
\right) \\
+ \int_{(\PP^1)^d} \log \left( \max_i \{ |f_i| \} \right) c_1(\overline{H}_1)
\wedge \cdots \wedge c_1(\overline{H}_d)
\end{multline*}
Hence, there is a constant $M_2$ independent on $P$ such that
\[
 \int_{(\PP^1)^d} \log (|f_i|) 
c_1(\overline{H}_1)
\wedge \cdots \wedge c_1(\overline{H}_d) \leq M_2
\]
for all $i$.
Thus, by \cite[Lemma~4.1]{MoArht},
we have our claim.

\bigskip
Let us consider a general case. We use the notation in the definition
of the largeness of a polarization.
Clearly, we may assume that $X = \PP^n_K$ and $L = \OO_{\PP^n_K}(1)$.
Let $K'$ be the function field of $B'$, and
$\overline{B}'$ a polarization of $K'$ given
by $(B'; \mu^*(H_1), \ldots, \mu^*(H_d))$.
Then, for all $x \in \PP^n(\overline{K})$,
\[
h^{\overline{B}'}_{L_{K'}}(x) = \frac{1}{[K':K]}h^{\overline{B}}_L(x).
\]
Thus, we may also assume that $B' = B$.
Moreover, there is a positive integer
$b$ with $\overline{L}_i^{\otimes b} \succsim \overline{\OO_{\PP^1_{\ZZ}}(1)}$
for every $i$. 
Hence, $\overline{H}_i^{\otimes b}
\succsim \nu^*(p_i^*(\overline{\OO_{\PP^1_{\ZZ}}(1)}))$.
Let 
$\overline{B}'$ be a polarization of $K$ given by
\[
\overline{B}' = \left(B; \nu^*(p_1^*(\overline{\OO_{\PP^1_{\ZZ}}(1)})), \ldots,
\nu^*(p_d^*(\overline{\OO_{\PP^1_{\ZZ}}(1)})) \right).
\]
Let $\overline{\OO_{\PP^n_B}(1)}$ be a $C^{\infty}$-hermitian line bundle
on $\PP^n_B$ given the pull-back of 
$(\OO_{\PP^n_{\ZZ}}(1), \Vert\cdot\Vert_{FS})$
via $\PP^n_B \to \PP^n_{\ZZ}$.
Then, since $\overline{\OO_{\PP^1_{\ZZ}}(1)}$ and $\overline{\OO_{\PP^n_B}(1)}$
are nef, 
we can see that, for all $x \in \PP^n(\overline{K})$,
\[
h^{\overline{B}'}_{(\PP^n_B, \overline{\OO_{\PP^n_{B}}(1)})}(x) \leq
b^d h^{\overline{B}}_{(\PP^n_B, \overline{\OO_{\PP^n_{B}}(1)})}(x).
\]
Thus, by the previous claim, we get our theorem.
\QED

\begin{Remark}
\label{rem:north:not:hold}
In order to guarantee Northcott's theorem,
the largeness of a polarization is crucial.
The following example shows us that even if the polarization
is ample in the geometric sense, Northcott's theorem does not hold.

Let $k = \QQ(\sqrt{29})$, $\epsilon = (5 + \sqrt{29})/2$, and
$O_k = \ZZ[\epsilon]$. We set
\[
E = \Proj\left( O_k[X, Y, Z]/(Y^2Z + XYZ + \epsilon^2YZ^2 - X^3) \right).
\] 
Then, $E$ is an abelian scheme over $O_k$. Then, as in the proof of
\cite[Proposition~3.1.1]{MoArht},
we can construct a nef $C^{\infty}$-hermitian line bundle $\overline{H}$ on $E$
such that $[2]^*(\overline{H}) = \overline{H}^{\otimes 4}$ and
$H_k$ is ample on $E_k$, $c_1(\overline{H})$ is positive on $E(\CC)$, and that
$\adeg\left(\acherncl_1(\overline{H})^2\right) = 0$.
Let $K$ be the function field of $E$. Then,
$\overline{B} = (E; \overline{H})$ is a polarization of $K$.
Here we claim that Northcott's theorem dose not hold 
for the polarization $(E, \overline{H})$ of $K$.

Let $p_i : E \times_{O_k} E \to E$ be the projection to the $i$-th factor.
Then, considering $p_2 : E \times_{O_k} E \to E$,
$(E \times_{O_k} E, p_1^*(\overline{H}))$ gives rise to
a model of $(E_K, H_K)$.
Let $\Gamma_n$ be the graph of $[2]^n : E \to E$, i.e.,
$\Gamma_n = \{ ([2]^n(x), x) \mid x \in E \}$.
Moreover, let $x_n$ be a $K$-valued point of $E_K$ arising from $\Gamma_n$.
Then, if we denote the section $E \to \Gamma_n$ by $s_n$, then
\begin{align*}
h_{H_K}^{\overline{B}}(x_n) & = \adeg 
\left( p_1^*(\overline{H}) \cdot p_2^*(\overline{H}) \cdot \Gamma_n \right)
= \adeg \left( s_n^*( p_1^*(\overline{H})) \cdot 
s_n^*(p_2^*(\overline{H})) \right) \\
& = \adeg \left( ([2]^n)^*(\overline{H}) \cdot \overline{H} \right)
= \adeg \left( \overline{H}^{\otimes 4^n} \cdot \overline{H} \right)
= 4^n \adeg \left(\overline{H} \cdot \overline{H} \right) = 0.
\end{align*}
On the other hand, $x_n$'s are distinct points in $E_K(K)$.
\end{Remark}

\end{document}